
\input amstex
\documentstyle{amsppt}
\magnification=\magstep1
\vsize 21 true cm
\hsize 15.0 true cm
\abovedisplayskip=1in
\belowdisplayskip=1in
\hoffset=0.05in
\voffset=0.05in

\catcode `\@=11
\let \logo@=\relax
\catcode  `\@=\active

\parindent=5mm
\NoRunningHeads
\NoBlackBoxes

\topmatter
\title
Castelnuovo Regularity  for  Smooth  Subvarieties of dimensions 3 and 4
\endtitle
\author
Sijong Kwak$^*$\endauthor 

\thanks $*$ Partially supported by  GARC(Global Analysis Research Center) 
and KIAS(Korea Institute for Advanced Study), Seoul, Korea.
\endthanks

\affil
Department  of  Mathematics \\ Columbia  University 
\endaffil


\endtopmatter
\document
\baselineskip=12pt
\subhead{\S 0. Introduction}\endsubhead

Let  $X$  be  a  smooth,  non--degenerate  projective  
variety  of  dimension  $n$  
and of degree  $d$  in  ${\Bbb P}^r$.  
We  say  that  $X$  is  {\it $k$--normal}  if  the  homomorphism
$$
H^0  ({\Bbb  P}^r , {\Cal  O}_{\Bbb  P}(k))  \rightarrow  
H^0 (X,  {\Cal  O}_X(k))
$$
is  surjective,  i.e.,  hypersurfaces  of  degree $k$  cut  out  a  complete  linear
system on $X$.  We  say  that  $X$  is  {\it $k$--regular}  if  
$H^i ({\Bbb  P}^r, {\Cal  I}_X(k-i))  =  0  \quad  \text{for  all}  \,\,  i  \ge  1$,
where  ${\Cal  I}_X$  is  the  sheaf  of  ideals  of  $X$ in ${\Cal  O}_{{\Bbb  P}^r}$.
It  is  easy  to  see  that  $X$  is  $(k+1)$-- regular  if  and  only  if 
$X$  is  $k$--normal  and  
$H^i (X , {\Cal  O}_X(k-i))  =  0  \quad  \text{for  all }  i   >   0$. 
Let  reg$(X)  =  \min  \{ k  \in  {\Bbb   Z}  \:   X  \,\,  \text{is}  \,\,  k\text{--regular} \}.$

The importance of ${k}$--regularity stems from the following
well-known results (\cite{Mu1}, lecture 14):
if  $X$  is  $k$--regular,  then  
the  saturated  ideal  $I_X$  is  generated  by  homogeneous  
polynomials  of  degree  at  most  $k$  and  hence  there  is    
no  $(k+1)$-secant  line  to  $X$.  Furthermore,  the  Hilbert  
polynomial  and  the  Hilbert  function  of $X$  have  the  same  
values  for  all  integers  $m  \ge  k-1$.

There is a well-known conjecture concerning the 
$k$--normality and $k$--regularity of $X$:

\proclaim{Regularity conjecture \cite{EG}, \cite{GLP}}
\roster
\item  $X$  is  $k$--normal  
$\text{for  all}  \,\,  k \ge \text{deg}(X)-\text{codim}(X)$.
\item $X$ is $k$--regular 
$\text{for  all}  \,\,  k  \ge  \text{deg}(X)-\text{codim}(X)  +  1$, \newline  
i.e.,  $\text{reg}(X)  \le  d- (r-n)  +  1$.   
\endroster
\endproclaim

As F. Zak has pointed out, the generalized Bertini theorem shows that a dimension $n$, degree $d$,
 non-degenerate variety $X$ in ${\Bbb P}^r$ can have at most a  $(d-(r-n)+1)$ 
secant line, thus yielding further creedence to the conjecture. 
The conjecture is also sharp 
in the sense that for all $n$, there exist varieties of dimension $n$
with regularity $(d-(r-n)+1)$.
 
It would already be interesting to establish the weaker bound 
$
\text{reg}(X)  \le  \text{deg}(X) 
$
for an arbitrary  dimensional, smooth variety $X$,
 since by projection it is possible to 
produce many equations of degree $d$ vanishing on $X$ -- cones of projections -- 
in fact, enough to cut out $X$
 locally ideal--theoretically.  This is a classical result, which is the  foundation
of the theory of Chow coordinates.  For a modern proof, see \cite{Mu2}. 
It is worth noting, however,
 that unpublished examples of Bayer-Pinkham show that the 
ideal generated by the cones of projection
can be strictly smaller, in degree d, than the saturated ideal,
 e.g.,  for elliptic curves.

The purpose of this paper is to give new bounds for regularity in dimensions 3 and 4
 which are only slightly worse than the optimal ones suggested by the conjecture. 
Our method will yield new
bounds up to dimension 14, but as they get progressively worse as the dimension goes up,
we have not written them down here. 

\proclaim{Theorem} Let $X$ be a smooth variety of dimension n and  of degree d in ${\Bbb P}^r$
\roster
\item"(i)"
If $n = 3$, then $\text{reg}(X) \le d-\text{codim}(X) + 2$. 
\item"(ii)"
If $n = 4$, then $\text{reg}(X) \le d-\text{codim}(X) + 5$.
\endroster 
\endproclaim

The main body of the proof is given in \S 3. Note that for $n = 3$ (resp. 4)
we are off by 1 (resp. 4) from the conjectured bound. In the case $n=4$,
our bound is precisely $d$ in the ``standard" case $r = 2n+1$.    

Here is a sketch of the history of the regularity conjecture.
 
0.1  Curves:  The  case of curves  in  $\Bbb P^3$  was  settled  by 
Castelnuovo  in
1893  (\cite {C}); that of curves  in  $\Bbb P^r, r\ge 3$, by  Gruson,  Lazarsfeld, 
and  Peskine in 1983 (\cite{GLP}).  

0.2   Surfaces: The  first bound for smooth surfaces 
was obtained by Pinkham(\cite{P})
using  geometric  methods. 
Shortly afterward Lazarsfeld(\cite{L})  obtained the full conjecture
using more cohomological techniques. 
It should be noted that Pinkham's proof relied on the classically
established ``fact" that the 3 inverse images above a triple point of a generic projection 
of a smooth surface in $\Bbb P^5$ to $\Bbb P^3$ lie in general linear position in their 
fiber (which is a $\Bbb P^2$).  
Vladimir Greenberg (\cite {L}, p.425) noticed
that the classic proofs of this fact
rely on an unproven general position assertion.  
More recently, it has become clear that the required ``fact" is  false
(despite the recent claim to the contrary by Ran \cite{Ran1}, 2.3).  
This issue will be discussed in detail in the forthcoming Columbia Ph.D.
thesis of Dobler (\cite {D}).   

Much of  Lazarsfeld's
argument for surfaces generalizes immediately to arbitrary dimension.
The only missing ingredient is 
information  on   the  length  and position of  the  fibers  of a general linear projection 
${\pi}_{\Lambda}: X^n  \rightarrow  {\Bbb  P}^{n+1}$  
in order to  separate  points  in  the  fiber  by  homogeneous
polynomials  of  some  given  degree.
Two results give partial results in this direction: 
first, the theorem of Mather (\cite {Ma1})
showing that for $n \le 14 $, generic linear projections of $X$ to 
$\Bbb P^{n+1}$ are {\it stable} enables us, in small dimension, to use
his results on stable mappings(\cite {Ma2}).  
Note however that Mather's results give no information on the position of
the points in the fiber, e.g., are they in general linear position? 
 A partial answer 
(in all dimensions) to the position question is given by Ran's 
(dimension $+2$)-secant lemma (\cite{Ran2}): the family  of  all
$(n+2)$-secant  lines  to a  $n$-dimensional  smooth  projective
variety   $X$  is  at  most $(n+1)$-dimensional. 

0.3  Threefolds:
In 1987, V. Greenberg(\cite{G}) modified Lazarsfeld's argument
to give the same bound 
as ours for threefolds.  
Unfortunately, at one step in the proof, 
there is an unsubstantiated claim which leaves 
the impression that the proof is incomplete.
In this note we follow Greenberg's proof closely, and justify the unsubstantiated claim. 
It should be noted that the bound we obtain in dimension 4 is stronger than Greenberg's.
   
In 1989, Z. Ran(\cite{Ran1})  claimed  the regularity conjecture for  smooth  threefolds  $X$  in   
${\Bbb  P}^r,  r \ge 9$.   
The proof relies on  differential geometric techniques to show that 
the  family  of  all  $4$-secant  lines  to  $X$  is  at  most  $4$-dimensional. 
 The method of proof definitely fails for embeddings in ${\Bbb  P}^7$ (the most important case)
and  ${\Bbb  P}^8$, 
which our method covers. For a further discussion of Ran's result, see \cite{D}. 
 
0.4  Above dimension 3 only weaker bounds were known:  
Mumford  showed  in 1984
that  reg$(X)  \le (n+1)(d-1)-n+1$,  where  dim$(X)=n$ and deg$(X)=d$  
(\cite{BM}) which  is   improved  to
$
\text{reg}(X)  \le  \min \{ e,n+1 \} (d-1) - n+1, \,e=\text{codim}X
$
in  \cite{BEL}. 

0.5 For a discussion of what is known in the small 
codimension cases that we do not treat here, we
refer to \cite{Kw}. 

\bigskip
{\bf  Acknowledgements.}   I  would  like  to  express  my  gratitude  to  
Professor  Henry  Pinkham  for  his  encouragement  and  stimulating  discussion  
with  geometric  insight  throughout  my  graduate  studies  at  Columbia  University  
where  this  work  was  prepared.  
In  addition,  I  thank  him  for  his  organization  of  
our  weekly  seminar  group  for  many  years and for editorial, historical and 
mathematical help in preparing this paper. 
As already noted, this paper 
also owes a great deal to
Greenberg's unpublished Ph.D. thesis (\cite{G}) written under the direction of Professor Pinkham,
and of course to Lazarsfeld's important paper(\cite{L}).

\subhead{\S 1. Basic background}\endsubhead
In  this  section,  we  recall  the  definitions and well-known 
results  which
are  used  in  the  following  sections.

\proclaim{Definition 1.1} For  a  coherent  sheaf  $\Cal F$ 
on ${\Bbb  P}_k^n$,  $\Cal F$  is  $m$--regular  if   
$\,\,  H^i ({\Bbb P}^n ,  {\Cal  F}(m-i))  =  0  \quad  \text{for   all}  \,\,  i > 0$
$\,\,$  and    reg$({\Cal  F})$  \,\,  is  defined  by  
$\inf  \{ m   \in   {\Bbb  Z}  \:  {\Cal  F}  \,\,  \text{is}  \,\,  m\text{--regular} \}.$
\endproclaim

\proclaim{Proposition 1.2} 
Let  $\Cal F$  be  a  coherent  sheaf  on  ${\Bbb  P}^n$.  
If  $\Cal  F$  is  $m$--regular,  then
\roster 
\item"(a)"  ${\Cal  F}(m)$  is  generated  by  its  sections.
\item"(b)"  The  multiplication  map  
$
H^0({\Cal  F}(m)) \otimes 
H^0 ({\Cal  O}_{\Bbb  P^n}(t)) \rightarrow  
H^0({\Cal  F}(m+t))
$  
is   surjective  for  $t \ge  0$.
\item"(c)"  $\Cal F$  is  $(m+t)$--regular  for  $t  \ge  0$.
\endroster
\endproclaim
\demo{Proof} See \cite{Mu1}, pg. 100.\qed
\enddemo
\proclaim{Proposition 1.3} 
\roster
\item"(a)"   
Let $\Cal F$  be  a  $p$--regular  vector  bundle  
  and  $\Cal G$  be  a  $q$--regular  vector  bundle  on  
${\Bbb  P}^n$ which is defined over an algebraically closed field of characteristic zero. Then
${\Cal  F}  \otimes  {\Cal  G}$  is  $(p+q)$ regular  and  $S^k({\Cal  F})$ 
and 
${\Lambda}^k({\Cal F})$  are  $(kp)$--regular.
\item"(b)"
 Let  $\Cal F$  be  a  coherent  sheaf  on  $\Bbb P^n$  and  
$\cdots  \rightarrow  
{\Cal  F}_i  \rightarrow  
\cdots  \rightarrow  
{\Cal  F}_0  \rightarrow  
{\Cal  F}  \rightarrow  0$  
be  an  exact  sequence  of  coherent  sheaves  on  ${\Bbb  P}^n$ 
such  that
${\Cal  F}_i$ is $(p+i)$--regular. Then $\Cal F$ is $p$--regular.
\endroster
\endproclaim
\demo{Proof} See \cite{L}, pg. 428.\qed
\enddemo

\proclaim{Theorem  1.4} 
Let  $X  \subset  {\Bbb  P}^r$  be  a  smooth  $n$-dimensional  subvariety  
and ${\pi}_{\Lambda} \,\, : X^n  \rightarrow   \overline{X} \subset  {\Bbb P}^{n+1}$
  a generic linear projection to a hypersurface. 
Let $\overline{X}_k$ be the closed set of points 
$ \overline{x} \in \overline{X}$ such that the scheme-theoretic length of the fiber 
${\pi}_{\Lambda}^{-1}( \overline{x})$ is at least k, and  $X_k$ the inverse image of 
$\overline{X}_k$, so $ X_{k+1} \subset X_{k}$ for all $k$.
Assume that $n \le 14$, so that we are in Mather's  ``nice" range.
Then $ X_{n+2}$ is empty, and $ X_{k}$ has dimension at most $n+1 - k$. If $X_{k}$ has the expected 
dimension $n+1 - k$, then outside of a closed subset of smaller dimension each fiber 
above $\overline{X}_k$ is a set of $k$ reduced points.
\endproclaim

\demo{Proof}
  This follows from the main theorem of \cite{Ma1} and the discussion in \S 5 of \cite{Ma2}.
A key ingredient is the inequality, for any $\overline{x} \in \overline{X}$
$$\sum_{x \in \pi^{-1}(\overline{x})} (\delta_x + \gamma_x) \le n+1$$
where $\delta_x$ is the length of the fiber at $x$ and $\gamma_x$ is another 
non--negative invariant defined (and 
computed) for all stable germs in the ``nice" range in \cite{Ma2}.
It is also known that $\delta_x = m_x$, 
where $m_x$ is the maximum number of points in the fiber
of the germ map in the neighborhood of $x$  (\cite {DG}). \qed
\enddemo  

In the classical literature (for example in Bertini's book \cite{B},
chapter 9, no.9, page 195 or more recently in \cite{Llu}) this result
is claimed {\it in all dimensions}.  The proof is incorrect, as
pointed out in \cite{G}, pg. 6 : assuming that $X$  is
$n$--dimensional and non--degenerate in ${\Bbb P}^{2n+1}$, the authors
implicitly assume  that the family of $n$--dimensional linear
subspaces of ${\Bbb  P}^{2n+1}$ intersecting $X$ in  at least $n+1$
points is an irreducible variety whose generic point is a linear span
of  $n+1$ points of $X$ in general position. It is easy to give
counterexamples to this assertion:  the following one, due to
Lazarsfeld, appears in \cite{G} :  let $E$ be a plane elliptic curve,
and embed $X = {\Bbb P}^1 \times E \subset {\Bbb P}^1 \times {\Bbb
P}^2$ by the Segre embedding  to ${\Bbb P}^5$.   Then the family of
trisecant planes has at least two components.  The generic
element of the first is the linear span of 3 points of $X$ in general position, while the generic
element of the second is the generic plane passing through 
the generic line in the Segre plane $x_0 \times {\Bbb P}^2$. 
(After this paper was prepared, Joel Roberts kindly brought to our
attention his unpublished 1969 Harvard Ph. D thesis, "Ordinary
Singularities of Projective Varieties", where the problem with the
classical method of proof is discussed in section 3.5.) 

It would be very interesting to know what actually happens outside the 
``nice" range of Mather.  See \cite{MP}, pg. 111 for a related
conjecture. In this direction we have a result of Ran, which we quote
but will not use in this paper:

\proclaim{Theorem  1.5 (The  dimension+2-secant lemma)}
Let  $X  \subset  {\Bbb  P}^r$  be  a  smooth  $n$-dimensional  subvariety  and  let  
$Y$  be  an  irreducible  variety  parametrizing  a  family of lines  in  
${\Bbb  P}^r$.  Assume  that,  for  a  general  $L_y$,  the  length  of  the  
scheme-theoretic  intersection  $L_y  \cap  X$  is  at  least   $n+2$.  Then  we  have  
$$
\text{dim}(\cup_{y  \in  Y}  L_y)  \le   n+1
$$
\endproclaim
\demo{Proof}
See  \cite{Ran2}.\qed
\enddemo

The importance of the result is that it is true in all dimensions. 
In the nice range it
is an immediate corollary of Theorem~1.4.

In this note, where we only treat the case $n \le 4$,
 we do not use the full strength of Theorem~1.4.  
Indeed  if $n \le 5$,  
all fibers of  ${\pi}_{\Lambda} \, : X^n  \rightarrow  
 \overline{X} \subset  {\Bbb P}^{n+1}$ 
as above are curvilinear, by an easy dimension count.  
This implies that
the result we need follows from \cite{PR}, \cite{MP} or Theorem~1.5 above.
Note that ``curvilinear" simply means that the finite
algebra (the algebra $Q(f_S)$ associated to the multigerm, defined in 
\cite{Ma2}, pg.188) appearing in the fiber of the projection has rank
at most 1,  in the language of \cite{Ma2}. 
In other words, it is a product of rings of the form 
${\Bbb C}[x]/(x^k)$. 
Some authors (\cite {MM}, \cite{MP})
say that  ${\pi}_{\Lambda}$ is {\it of corank $1$}.

\subhead{\S 2. A Construction of Lazarsfeld \cite{L}}
\endsubhead

We first establish some notation that we will use throughout the rest of this note.

Let $X$ be a smooth non--degenerate complex projective variety of degree $d$ 
and of dimension $n$ in ${\Bbb P}^r$. 
Let $\Lambda^{r-n-2}$ be a general linear space of
dimension $(r-n-2)$, disjoint from $X$, and let $p_1:Bl_{\Lambda}{\Bbb P}^r  \longrightarrow
 {\Bbb P}^r$ be the blow--up of 
${\Bbb  P}^r$ along $\Lambda$. 
Let $p_2:Bl_{\Lambda}{\Bbb P}^r  \longrightarrow
 {\Bbb P}^{n+1}$ be the projection associated to $\Lambda$. 
We may realize the projection by choosing any linear subspace $\Psi \subset {\Bbb P}^r$ 
of dimension $n+1$ not intersecting $\Lambda$. 
In the intrinsic realization of a projection,
 a point in the projected space corresponds to a linear subspace $H$ of ${\Bbb P}^r$, of dimension
 $r-n-1$, containing $\Lambda$.
 After choosing a $\Psi$, we can parametrize 
the projection by the (unique)
point of intersection $q = \Psi \cap H$. 
Going the other way, we can reconstruct
$H$ from $q \in \Psi$ by taking $H = \langle\Lambda,q\rangle$.  
Obviously the projection does not depend on the choice of $\Psi$. 
In the next 
section we will vary $\Psi$ without changing $\Lambda$, in order to ensure
that $\Psi$ does not vanish at a given finite set of points. 
This is obviously 
always possible, and amounts to a change of coordinates. We choose homogeneous 
coordinates $T_0,T_1\dots,T_r$ on ${\Bbb P}^r$ so that 
$\Lambda = Z(T_0,T_1\dots,T_{n+1})$, the locus where $T_0,T_1\dots,T_{n+1}$ all vanish, and
$\Psi = Z(T_{n+2}, \dots,T_r)$.  
Setting $V={\Bbb C}\cdot T_{n+2}\oplus {\Bbb C}\cdot T_{n+3}\oplus\dots
\oplus {\Bbb C}\cdot T_{r}$, and 
$W={\Bbb C}\cdot T_{0}\oplus {\Bbb C}\cdot T_{1}\oplus\dots\oplus {\Bbb C}\cdot T_{n+1}$,  
we also write 
$\Lambda = {\Bbb P} (V)$ and  $\Psi = {\Bbb P} (W)$. 
Note that
$
Bl_{\Lambda}{\Bbb P}^r={\Bbb  P}({\Cal O}_{\Psi}(1)\oplus
(V\otimes {\Cal O}_{\Psi}))=\{(x,q):\,x \in L_q=\langle\Lambda,q\rangle\,\,, q \in \Psi\}
$. 

Consider the   diagram
$$
\CD
Bl_{\Lambda}{\Bbb  P}^r  =  
\Bbb  P({\Cal  O}_{\Psi}(1)  \oplus( V\otimes \Cal  O_{\Psi})) @>{p_2}>> 
\Psi = \Bbb  P^{n+1} \\  
     @V{p_1}VV  @.  \\
X  \subset  \Bbb  P^r   @. 
\endCD
\tag2.1
$$

\noindent Let $\pi_{\Lambda}$ be the restriction of  
$p_2  \circ  p^{-1}_1$ to  $X$. 
By the choice of $\Lambda$, $\pi_{\Lambda}$ is
a (finite) morphism.

From the exact sequence
$0 \rightarrow {\Cal I}_{X/{\Bbb P}^r}(k)\rightarrow {\Cal O}_{{\Bbb P}^r}(k)
\rightarrow {\Cal O}_X(k)\rightarrow 0$
for each $k\in\Bbb Z$,  we have the following sequence
$$
0\rightarrow
p_{2\ast}(p_1^{\ast}\Cal I_{X/{\Bbb P}^r}(
k)) \rightarrow  p_{2\ast}(p_1^{\ast}\Cal O_{{\Bbb P}^r}(k))
\overset{\omega_{n,k}}\to 
\rightarrow p_{2\ast}(p_1^{\ast}{\Cal O_X}(k))
\rightarrow
R^1p_{2\ast}(p_1^{\ast}{\Cal I}_{X/{\Bbb P}^r}(k))\rightarrow 0$$
because $R^1p_{2\ast}(p_1^{\ast}{\Cal O}_{{\Bbb P}^r}(k))=R^1p_{2\ast}({\Cal O}_{{\Bbb P}(\Cal E)}(k))
=0$, where $\Cal E={\Cal O}_{\Psi}(1)\oplus(V\otimes{\Cal O}_{\Psi}).$

Note that $p_{2\ast}(p_1^{\ast}{\Cal O_X}(k))={\pi}_{\Lambda\ast}{\Cal O}_X(k)$ and
$p_1^{\ast}{\Cal O}_{{\Bbb  P}^r}(1)={\Cal O}_{{\Bbb P}(\Cal E)}(1)$ which is the tautological line 
bundle.

The main issue is to prove the surjectivity of $\omega_{n,k}$ for some $n$
and $k$, where $\omega_{n,k}$ is the map in the above exact sequence.
By Nakayama's lemma, it suffices to show that for all $q \in \Psi$,

$$
\CD
p_{2\ast}(p_1^{\ast}{\Cal O}_{{\Bbb  P}^r}(k))\otimes {\Bbb C}(q) 
@>{\omega_{n.k}\otimes {\Bbb C}(q)}>> p_{2\ast}(p_1^{\ast}\Cal O_X(k))  
\otimes {\Bbb C}(q) \\
 @V{\cong}VV  @V\cong VV\\ 
 H^0(L_q,{\Cal O}_{L_q}(k))    @>>>    
 H^0(L_q,{\Cal O}_ {\pi^{-1}_{\Lambda}(q)} (k))  
 \endCD
\tag2.2
$$
is surjective. Equivalently, using the bottom horizontal arrow,
 it is enough to show that the finite scheme $\pi^{-1}_{\Lambda}(q)$ in 
$L_q=\langle\Lambda,q\rangle$ is $k$--normal for all $q\in \Psi$. 
Therefore,
we need some information on the finite schemes $\pi^{-1}_{\Lambda}(q)$ appearing in
the fiber of a generic projection $\pi_{\Lambda}: X \rightarrow  \Psi$.
 Note that a finite scheme of length $k$ is $(k-1)$--normal and it fails to be $(k-2)$--normal if and only if 
it is contained in a line.  
When the fiber is {\it reduced},
 to prove surjectivity in degree $k$
it is enough, for each point $p$ in the fiber, 
to find a polynomial of degree $k$ vanishing at all the points
in the fiber
other than $p$, and not vanishing at $p$. 

Suppose  that  $\omega_{n,k}$ is surjective  for  some  $k>0$.  
Recall  that   
$$
\align
&p_{2\ast}(p_1^{\ast}{\Cal  O}_{{\Bbb  P}^r}(k))=Sym^k({\Cal
O}_{\Psi}(1)  \oplus ( V  \otimes  {\Cal  O}_{\Psi}))\\   =&{\Cal
O}_{\Psi}(k)  \oplus  ( V  \otimes  {\Cal  O}_{\Psi}(k-1))  \oplus  
(S^2(V)  \otimes  {\Cal  O}_{\Psi}(k-2))  \oplus  
\cdots  \oplus  
(S^k(V)  \otimes  {\Cal  O}_{\Psi}),
\endalign
$$
where  $S^i(V)$  is  the  $i$-th  symmetric  power  of $V$.  
After  twisting  by  $(-k)$,  we  get  the  exact  sequence  
$$
\CD
0  \rightarrow  E_{n,k}  \rightarrow  
(S^k(V)  \otimes  {\Cal  O}_{\Psi}(-k)) \oplus  \cdots  \oplus  
(V  \otimes  {\Cal  O}_{\Psi}(-1)) \oplus  
{\Cal   O}_{\Psi}  @>  \tilde{\omega}_{n,k} >>  
\pi_{\Lambda_{\ast}}{\Cal  O}_X  \rightarrow  0, 
\endCD
$$  
where  $E_{n,k}  =  Ker ( \tilde{\omega}_{n,k})$.  
Since  $\pi_{\Lambda_{\ast}}{\Cal  O}_X$  is  a locally  
Cohen-Macaulay module of codimension 1
for ${\Cal  O}_{\Psi}$,  $E_{n,k}$  is  locally  free.
 
Now, following Greenberg (\cite{G}), we generalize this construction. 
For each $j \ge 1$ we
choose a subspace $V_j$ of $S^j(V)$ such that the map $\tilde{\omega}_{n,k, {\Cal V}}$ 
obtained by restriction,
$$
\CD
(V_k  \otimes  {\Cal   O}_{\Psi}(-k)) \oplus  \cdots   \oplus 
(V_1  \otimes  {\Cal  O}_{\Psi}(-1))  \oplus  
{\Cal   O}_{\Psi}  @>  \tilde{\omega}_{n,k,{\Cal V}} >>  
\pi_{\Lambda_{\ast}}{\Cal  O}_X 
\endCD
\tag2.3 
$$
is still surjective. 
We write  ${\Cal V}$ for the collection of $V_j$. 
We will always take
$V_1 = V$ and $V_2 = S^2(V)$.  
As before, the kernel $E_{n,k, {\Cal V}}$ of $\tilde{\omega}_{n,k,{\Cal V}}$ is 
locally free.
  
\proclaim{Lemma  2.1} If $\tilde{\omega}_{n,k, {\Cal V}}$ is surjective, then:  
\roster
\item  
If  $H^1({E}_{n,k, {\Cal V}}(m))  =  0$,  then  $X$  is  $m$-normal.
\item  
$reg(X)  \leq  reg(E_{n,k,{\Cal V}})  =  reg(\Lambda^{\text{rank}(E_{n,k,{\Cal V}})-1}(E_{n,k,{\Cal V}}^{\ast})  
\otimes  detE_{n,k,{\Cal V}})$\newline 
$\leq  (\text{rank}(E_{n,k,{\Cal V}})-1)  reg(E_{n,k,{\Cal V}}^{\ast})  -  c_1(E_{n,k,{\Cal V}})$. 
\endroster
\endproclaim
\demo{Proof} $(1)$ is a straightforward extension of  the proof of
 lemma 1.5 of \cite{L}.
For $(2)$, first check  that reg$(X) \leq {\text{reg}(E_{n,k,{\Cal V}})}$. 
Suppose  that  $E_{n,k,{\Cal V}}$  is  $(m+1)$--regular, 
in other  words,   
$h^i(\Psi,  E_{n,k,{\Cal V}}(m+1-i))  =  0,  \, i>0$.  
We can  show  that  $X$  is  $(m+1)$--regular,  equivalently,  $X$  is  $m$--normal  and  
$h^i({\Bbb  P}^r,  {\Cal  O}_X(m-i))  =  0,  \,\,  i>0$.  
From the exact sequence induced from (2.3),
$$
\CD
0\rightarrow
E_{n,k, {\Cal V}}\rightarrow
V_k  \otimes  {\Cal   O}_{\Psi}(-k)  \oplus   
\cdots   \oplus 
V_1  \otimes  {\Cal  O}_{\Psi}(-1)  \oplus  
{\Cal   O}_{\Psi}  @>  \tilde{\omega}_{n,k,{\Cal V}} >>  
\pi_{\Lambda_{\ast}}{\Cal  O}_X
\rightarrow
0  
\endCD
$$
it is easy to check that
$H^i(\Psi,  E_{n,k,{\Cal V}}(m+1-i))  =  0,  \,\,  i>0$  implies  $H^i({\Bbb  P}^r,$ 
${\Cal  O}_X(m-i)) = 0$  for  $i>0$.  
Therefore,  $X$  is  $(m+1)$--regular and  reg$(X)\leq$ reg$(E_{n,k,{\Cal V}})$. 
The second inequality follows from the isomorphism \newline
$E_{n,k,{\Cal V}}  \simeq  \Lambda^{\text{rank}(E_{n,k, {\Cal V}})-1}(E_{n,k,{\Cal V}}^{\ast})  
\otimes  detE_{n,k,{\Cal V}}$ and 
Proposition~1.3 (a).\qed 
\enddemo

\proclaim{Lemma 2.2}
reg$(E_{n,k,{\Cal V}}^{\ast})   \le  (-2)$
\endproclaim
\demo{Proof}  See \cite{L}, Lemma 2.1. where $A^{\ast}$ is
(--1)-regular and  $B^{\ast}$ is (--2)-regular. From the sequence
$
0  \rightarrow  
\tilde{A}^{\ast}  \rightarrow  
\tilde{B}^{\ast}  \rightarrow  
E_{n,k,{\Cal V}}^{\ast}  \rightarrow
0
$
in  the same lemma of \cite{L}
and Proposition~1.3 (b), $E_{n,k,{\Cal V}}^{\ast}$ is (--2)-regular.
\qed 
\enddemo

\proclaim{Proposition 2.3 (Greenberg)} Assume 
 $\tilde{\omega}_{n,k, {\Cal V}}$ is surjective.
Then 
$$
reg(X)  \leq \text{deg}(X)-\text{codim}(X)  +  1 + \sum_{j=3}^{k}  (j-2) \dim{V_j}
$$
\endproclaim
\demo{Proof}
This is an easy computation obtained from Lemma 2.1,(2).
\qed 
\enddemo

The goal is therefore to get the $V_j$ as small as possible. 
Curiously, the quadratic polynomials
$V_2$ do not contribute to the final result. 
To obtain the regularity
conjecture one must have $V_j = 0$ for all $j \ge 3$.  
This is what happens for 
smooth curves and surfaces, see \cite{L} pg. 425.

\subhead{\S 3. Castelnuovo regularity for smooth subvarieties of dimensions 3, 4}
\endsubhead

    We first deal with the case $n = 3$.  
By Mather's Theorem~1.4, we see that 
outside of the finite set of points $\overline{X}_4$ (using the notation of that theorem) the fibers of the
general projection have length at most three, and so can be separated by quadratic polynomials.
Above the points $\overline{X}_4$, the fiber consists of 4 reduced points and this too can be
 separated by polynomials of degree 2 {\it unless the four points are aligned}.  
In the next theorem, we
show that this can be handled by a subspace $V_3$ of dimension 1. 
  
\proclaim{Theorem 3.1}
Let $X$ be a smooth 3-fold of degree d in ${\Bbb P}^r$. Then we have 
$\text{reg}(X)  \le  \text{deg}(X)-\text{codim}(X)  +  2$  
\endproclaim
\demo{Proof}
We use the notation of $\S 2$, simply setting $n=3$.   
 
By Theorem 1.4,  the finite scheme $\pi^{-1}_{\Lambda}(q)$ in
$\L_q=\langle\Lambda,q\rangle$  has length at most 4 for all $q \in \Psi$. 
As in Theorem~1.4, let $\overline{X}_k =\{ q  \in \overline{X}\,|$
the length of $\pi^{-1}_{\Lambda}(q) \ge k\}$, and let $Y_4$ be the subset of $\overline{X}_4$
above which the fiber consists of 4 distinct, collinear points. 
Then the commutative diagram (2.2) is surjective for all
$q  \notin  Y_4 $ and $k=2$.  
In other words, every fiber $\pi^{-1}_{\Lambda}(q)$, $q \notin Y_4$
is  2--normal in  $\L_q=\langle\Lambda,q\rangle$. 
 Therefore,
the morphism  $\omega_{3,2}:{\Cal  O}_{\Psi}(2)  \oplus  
(V  \otimes  {\Cal  O}_{\Psi}(1))\oplus  
(S^2(V)  \otimes  {\Cal  O}_{\Psi}) \rightarrow  \pi_{\Lambda_{\ast}}{\Cal  O}_X (2)$  
is surjective for all $y  \notin Y_4$.
Thus the morphism 
$$\omega_{3,2}\otimes{\Cal  O}_{\Psi}(1):{\Cal  O}_{\Psi}(3)  \oplus  
(V  \otimes  {\Cal  O}_{\Psi}(2))\oplus  
(S^2(V)  \otimes  {\Cal  O}_{\Psi}(1))   \rightarrow  \pi_{\Lambda_{\ast}}{\Cal  O}_X (3)
$$ 
is  also surjective for all 
$y  \notin Y_4$.
 Recall that  $\overline{X}_4$, and therefore its subset $Y_4$  is 
(at most) a finite set by Theorem~1.4. 
Write  $Y_4 = \{q_1,\dots,q_t\}$, so
the points $\pi^{-1}_{\Lambda}(q_i)$
are collinear. 
Denote by $\ell_{q_i}$ the line supporting the fiber above $q_i$, i.e., 
$\pi^{-1}_{\Lambda}(q_i) \subset  \ell_{q_i} \subset   L_{q_i}$ for all $q_i \in Y_4$.  
Let $p_i$ be the point of intersection of the line  $\ell_{q_i}$ with $\Lambda$.
Now, choose a linear form $H(T_5, \dots, T_r)$ on $\Lambda$
that does not vanish on  $\{p_1,p_2,\dots p_t\}$. 
Viewing, in the obvious way, H as a form on $\Bbb P^r$,  
we can restrict $H$ to $\ell_{q_i}$ and by construction, 
$H$  vanishes only at the point $q_i$. 
Let $U(T_0, \dots, T_4)$ be a linear form on $\Psi$
that does not vanish on  $\{q_1,q_2,\dots q_t\}$. 
Then $U$ restricted to the line  $\ell_{q_i}$
only vanishes at $p_i$. 
Clearly $\{H,U\}$, when restricted to $\ell_{q_i}$, is a homogeneous coordinate system  for all $i=1,2,\dots,t$,  
and $\pi^{-1}_{\Lambda}(q_i) \subset  \ell_{q_i}$, 
for all $q_i \in Y_4$, can be separated by cubic polynomials in $\{H,U\}$.  
Therefore, the natural morphism
$$
\lbrack{\Cal  O}_{\Psi}(3)  \oplus  
(V  \otimes  {\Cal  O}_{\Psi}(2))\oplus  
(S^2(V)  \otimes  {\Cal  O}_{\Psi}(1)) \oplus 
(H^3 \otimes {\Cal  O}_{\Psi}) \rbrack\otimes {\Bbb C}(q)
 \rightarrow  \pi_{\Lambda_{\ast}}{\Cal  O}_X (3)\otimes {\Bbb C}(q)
$$ 
is surjective for all $q \in Y_4$
 because  the left-side of the morphism generates all cubic polynomials
in the variables $\{H,U\}$ on the line.  
In other words,
taking $V_3$ to be the one-dimensional space generated by $H^3$, 
$\omega_{3,3, {\Cal V}}$ is surjective so
by Proposition~2.3, the Theorem is proved. 
\qed
\enddemo

We now turn to the case of dimension 4.

\proclaim{Theorem 3.2}
Let $X$ be a smooth, non-degenerate subvariety of degree d and
of dimension $4$ in $\Bbb P^r$. Then we have
reg$(X) \le d-codimX+5$  
\endproclaim

We will show that we can take $U_3$ of dimension 2, $U_4$ of dimension 1, 
and all higher $U$ of dimension 0. 
The codimension 2 case is dealt with in \cite{Kw}, where we get 
reg$(X) \le d-1$ for smooth threefold or fourfold $X$  
and will not be discussed here, so we assume $n \ge 7$.

\demo{Proof} As in $\S 2$, take a generic projection 
${\pi}_{\Lambda} : X^4  \rightarrow   \overline{X} \subset  \Psi = {\Bbb P}^{5}$.

Let $\overline{X}_k$ and $X_k$ have the same meaning as in Theorem~1.4, and let 
$Y_4$ be the subset of $\overline{X}_4$ where the fiber contains 4 distinct {\it collinear} 
points. 
Let $Y_5 = Y_4 \cap \overline{X}_5$. 
By Theorem~1.4,
dim$Y_k  \le  5-k$. So,  $Y_4$  is (at most) a 1-dimensional subvariety in 
$\pi_{\Lambda}(X)  \subset \Psi$ 
and $Y_5 =\{q_1,\dots,q_t\}$ is a finite set of $t$ points. 
For each point $q \in Y_4$, there
is a unique line $\ell_q$ which is the support of the collinear points in
the  fiber above $q$. As $q$ varies in $Y_4$, 
the lines
$\ell_q$ intersect $\Lambda$ in a variety $Z$ which is at most 1--dimensional.
As $q_i$ varies in $Y_5$, the lines
$\ell_{q_i}$ intersect $\Lambda$ in points $p_i$. 
Choose  linear forms $H_1 ,H_2$ on $\Lambda$
such that $H_1$ is nonzero  at $\{p_1,p_2,\dots p_t\}$ and $H_2$ 
is nonzero at $H_1\cap\  Z =\{t_1,\dots,t_m\} \subset \Lambda$.
The morphism in the diagram (2.2) 
is surjective for all $q  \notin  Y_4$.  
 Therefore,
the morphism  $\omega_{4,2}:{\Cal  O}_{\Psi}(2)  \oplus  
V  \otimes  {\Cal  O}_{\Psi}(1)\oplus  
S^2(V)  \otimes  {\Cal  O}_{\Psi} \rightarrow  \pi_{\Lambda_{\ast}}{\Cal  O}_X (2)$  
is surjective for all $q  \notin  Y_4$.
Note then that the morphism 
$$\omega_{4,2}\otimes{\Cal  O}_{\Psi}(1):{\Cal  O}_{\Psi}(3)  \oplus  
(V  \otimes  {\Cal  O}_{\Psi}(2))\oplus  
(S^2(V) \otimes {\Cal  O}_{\Psi}(1))  \rightarrow  \pi_{\Lambda_{\ast}}{\Cal  O}_X (3)
$$ 
is  also surjective for all $q  \notin  Y_4$. 
We first treat the fibers above $Y_4 \setminus Y_5$. 
By construction, one of $H_1$ and $H_2$
acts as the second homogeneous coordinate on the support of the fiber, which is a line.  
So if we let $V_3 $ be the 2--dimensional subspace of $S^3(V)$ generated by $H_1^3$ and $H_2^3$, 
the corresponding map $\tilde{\omega}_{4,3, U}$ is surjective everywhere except possibly above $Y_5$.
 
Above $Y_5$, we have two kinds of fibers: 
\roster
\item 5 reduced collinear points (this is 4--normal)  
\item 4 reduced collinear points plus one extra point off the line in $L_{q_i}$. 
\endroster

 We note that there is no fiber which contains four aligned points, and one of the 
points in the fiber is nonreduced and has multiplicity 2. This case is ruled out by 
J.Mather's inequality that we quote on page 4(as in the case of surfaces).
 
In the first case, it is necessary to include the element $H_1^4$ in $V_4$
 to obtain surjectivity of the corresponding $\tilde{\omega}_{4,4, {\Cal V}}$. 
Note that
by construction $H_1$ only vanishes at the point $q_i$ on the line $\ell_{q_i}$.
The argument is the same as in dim$(X)$=3.

The second case is more delicate. 
To prove the Theorem,
 we must avoid adding any extra 
global sections to $V_3$ or $V_4$. 
Since we can work one fiber at a time, we drop the $i$ index.
The fiber above $q \in Y_5 \subset \Psi$ spans  
a 2-dimensional linear space in ${L_{q}}$\
 which meets $\Lambda$ in a line 
$N_{q}$ containing the point $p \in \Lambda$. 
Since we have only a finite number of
$N_{q}$, we may of course assume that neither $H_1$ nor $H_2$
 vanishes identically on $N_{q}$.
We may also assume that $q \notin \pi^{-1}_{\Lambda}(q)$, 
in other words that $q \notin X \cap \Psi$,
 by moving $\Psi$ as explained in $\S 2$, 
since $Y_5$ is a finite set.
Thus the fiber ${\Bbb P}^2$ has homogeneous coordinates 
$(U, H_1, H_2)$, where the line $U = 0$
is the intersection of the fiber with $\Lambda$. 
Thus we can assume that the four aligned points have coordinates
$(u_i, a,b)$ and the fifth point has coordinates $(u, c, d)$.
By construction, none of the $u_i$
or $u$ are equal to 0. 
We may choose the $H_j$ so that the same is true for $a, b, c, d$. 
 We also know that none of the points is the point $q = (1,0,0)$.
  
We need to prove the following easy lemma, 
where for simplicity of notation we set $H_1 = X$ and $H_2 = Y$.
\proclaim{Lemma 3.3}
Let $U, X, Y$ be homogeneous coordinates on ${\Bbb P}^2$, and suppose given 5 points, 
$p_i = (u_i, a,b)$, $1 \le i \le 4$ and $ p_5 = (u, c, d)$. 
 Assume that none of the $u_i, u, a, b, c, d$ are 0.  
Then the points can be separated 
using only cubics whose equation only contain the monomials 
$$U^3, U^2X, U^2Y, UX^2,UXY, UY^2 , X^3, Y^3,$$ 
\endproclaim
\demo{Proof} This is easy to check by elementary linear algebra. 
By symmetry it is enough
to construct a cubic that vanishes at $p_1$, $p_2$, $p_3$, and $p_5$. 
(Note that it is trivial to construct
a linear form vanishing at the four aligned points, but not at $p_5$.)
Consider the  cubic polynomial
$$
U^3+(a_{1,0}U^2X+a_{0,1}U^2Y)+(a_{2,0}UX^2+a_{1,1}UXY+a_{0,2}UY^2)+(a_{3,0}X^3+a_{0,3}Y^3)=0
$$
First, note that this cubic doesn't vanish identically on the line $\ell = \{aY - bX = 0\}$
 containing the 4-collinear points 
because it doesn't vanish at the point $(1,0,0)$ which is on $\ell$. 
Evaluating $X = a$, $Y = b$ ,  we must have the cubic polynomial in $U$:
$$
U^3+(a_{1,0}a+a_{0,1}b)U^2+(a_{2,0}a^2+a_{1,1}ab+a_{0,2}b^2)U^2+(a_{3,0}a^3+a_{0,3}b^3)
$$
This must be equal to $(U-u_1)(U-u_2)(U-u_3)$ if the polynomial is to vanish at $p_1$, $p_2$, $p_3$.
For this to happen, equating coefficients, 
we must solve a system of three linear equations in the seven
unknowns $a_{1,0}, a_{0,1},a_{2,0}, a_{1,1}, a_{0,2}, a_{3,0}, a_{0,3} $. 
This yields a $4$--dimensional
family of solutions. 
 Finally we must force the solution to pass through $p_5$. 
This is one extra linear condition.
So we have found a cubic (in fact a 3--dimensional family) passing though $p_1$, $p_2$, $p_3$ and $p_5$.
  Since it does not contain the line $\ell$,
and since we have all its intersections with $\ell$, 
it does not contain $p_4$. 
Therefore, the point $p_4$ can be separated from the others.
\qed
\enddemo

We have shown that taking $V_3 = \langle H_1^3, H_2^3\rangle $ of dimension 2, 
$V_4 = \langle H_2^4\rangle$ of dimension 1, 
$\tilde{\omega}_{4,4, {\Cal V}}$ is surjective,
and therefore, by Proposition~2.3, the theorem is proved.\qed
\enddemo

\bigskip\noindent
{\bf Current Address:}
{School of Mathematics, Korea Institute for Advanced Study, 
207-43 Chungryangri-dong,Dongdaemoon-gu, Seoul 130-010, Korea}\newline
{\bf E-mail:} {\it kwak$\@$math.snu.ac.kr }
\enddocument